\documentclass[10pt,onecolumn]{article}

\usepackage{times}
\usepackage{latexsym}

\usepackage{amsmath} 
\usepackage{amssymb} 

\makeatletter
\def\@maketitle{%
  \vbox to 6.5cm{%
    \hsize\textwidth
    \linewidth\hsize
    \vspace{1.5cm}
    \centering
    {\bfseries\LARGE \@title \par}
    \vspace{12pt}
    {\fontsize{11pt}{13pt}\selectfont \begin{tabular}[t]{c}\@author \end{tabular}\par}
    \vfill} 
}
\renewcommand\section{\@startsection{section}{1}{\z@}%
                       {-12\p@ \@plus -4\p@ \@minus -4\p@}%
                       {6\p@ \@plus 4\p@ \@minus 4\p@}%
                       {\normalfont\large\bfseries
                        \rightskip=\z@ \@plus 8em\pretolerance=10000 }}
\renewcommand\subsection{\@startsection{subsection}{2}{\z@}%
                       {-12\p@ \@plus -4\p@ \@minus -4\p@}%
                       {6\p@ \@plus 4\p@ \@minus 4\p@}%
                       {\normalfont\fontsize{11pt}{13pt}\selectfont\bfseries
                        \rightskip=\z@ \@plus 8em\pretolerance=10000 }}
\renewcommand\subsubsection{\@startsection{subsubsection}{3}{\z@}%
                       {-12\p@ \@plus -4\p@ \@minus -4\p@}%
                       {6\p@ \@plus 4\p@ \@minus 4\p@}%
                       {\normalfont\normalsize\itshape}}
\renewcommand\paragraph{\@startsection{paragraph}{4}{\z@}%
                       {-12\p@ \@plus -4\p@ \@minus -4\p@}%
                       {-0.5em \@plus -0.22em \@minus -0.1em}%
                       {\normalfont\normalsize\itshape}}
\makeatother

\renewenvironment{abstract}%
  {\small
    \list{}{\labelwidth0pt
      \leftmargin0pt \rightmargin\leftmargin
      \listparindent\parindent \itemindent0pt
      \parsep0pt
      }%
    \item[\hskip\labelsep\bfseries\abstractname\enspace --] \itshape}{\endlist}

\newcommand{\keywordsname}{Keywords}
\newenvironment{keywords}%
  {\small
    \list{}{\labelwidth0pt
      \leftmargin0pt \rightmargin\leftmargin
      \listparindent\parindent \itemindent0pt
      \parsep0pt
      }%
    \item[\hskip\labelsep\bfseries\keywordsname:]}{\endlist}

\begin{document}

\pagestyle{plain} 

\title{Infinite classes of counter-examples \\
to the Dempster's rule of combination}

\author{Jean Dezert\\
ONERA, Dtim/ied\\
29 Av. de la  Division Leclerc \\
92320 Ch\^{a}tillon, France.\\
Jean.Dezert@onera.fr\\
\and
Florentin Smarandache\\
Department of Mathematics\\
University of New Mexico\\
Gallup, NM 87301, U.S.A.\\
smarand@unm.edu}
\date{}

\maketitle

\begin{abstract}
This paper presents several classes of fusion problems which cannot be directly attacked by the classical mathematical theory of evidence, also known as the Dempster-Shafer Theory (DST) either because the Shafer's model for the frame of discernment is impossible to obtain or just because the Dempster's rule of combination fails to provide coherent results (or no result at all). We present and discuss the potentiality of the DSmT combined with its classical (or hybrid) rule of combination to attack these infinite classes of fusion problems.
\end{abstract}



\section{Introduction}

In this paper we focus our attention on the limits of the validity of the Dempster's rule of combination in the Dempster-Shafer theory (DST) \cite{Shafer_1976}. We provide several infinite classes of fusion problems where the Dempster rule of combination fails to provide coherent results and we show how these problems can be attacked directly by the recent theory of plausible and paradoxical reasoning developed by the authors in \cite{Dezert_2002b,Dezert_2003,Dezert_Smarandache_2003,Dezert_2003f,Smarandache_2002,Dezert_2004Book}. We don't present in details here these two theories but just remind in next section the main differences between them. DST and DSmT are based on a different approach for modeling the frame $\Theta$ of the problem (Shafer's model versus free-DSm, or hybrid-DSm model), on the choice of the space (classical power-set $2^\Theta$ versus hyper-power set $D^\Theta$) on which will be defined the basic belief assignment functions $m_i(.)$ to be combined, and on the fusion rules to apply (Dempster rule versus DSm rule or DSm hybrid rule of combination).

\section{A short DST and DSmT presentation}

The DST starts by assuming an exhaustive and exclusive frame of discernment of the problem under consideration $\Theta=\{\theta_1,\theta_2,\ldots,\theta_n\}$. This corresponds to the Shafer's model of the problem usually denoted $\mathcal{M}^0(\Theta)$ \cite{Dezert_2004Book}.  The Shafer's model assumes actually that an ultimate refinement of the problem is possible so that $\theta_i$ are well precisely defined/identified in such a way  that we are sure that they are exclusive and exhaustive. From this Shafer's model,  a basic belief assignment (bba) $m_i(.): 2^\Theta \rightarrow  [0, 1]$  such that $m_i(\emptyset)=0$ and $\sum_{A\in 2^\Theta} m_i(A) = 1$ associated to a given body of evidence $\mathcal{B}_i$ is defined, where $2^\Theta$ is called the {\it{power set}} of $\Theta$, i.e. the set of all subsets of $\Theta$. The set of all propositions $A\in 2^\Theta$ such that $m_i(A)>0$ is called the core of $m_i(.)$ and is denoted $\mathcal{K}(m_i)$. Within DST, the fusion (combination) of two independent sources of evidence $\mathcal{B}_1$ and $\mathcal{B}_2$ is obtained through the following Dempster's rule of combination  \cite{Shafer_1976} : $[m_{1}\oplus m_{2}](\emptyset)=0$ and $\forall B\neq\emptyset \in 2^\Theta$:
 \begin{equation}
 [m_{1}\oplus m_{2}](B) = 
\frac{\sum_{X\cap Y=B}m_{1}(X)m_{2}(Y)}{1-\sum_{X\cap Y=\emptyset} m_{1}(X) m_{2}(Y)} 
\label{eq:DSR}
 \end{equation}
 
The notation $\sum_{X\cap Y=B}$ represents the sum over all $X, Y \in 2^\Theta$ such that $X\cap Y=B$. 
The orthogonal sum $m (.)\triangleq [m_{1}\oplus m_{2}](.)$ is considered as a basic belief assignment if and only if the denominator in equation \eqref{eq:DSR} is non-zero. The term $k_{12}\triangleq \sum_{X\cap Y=\emptyset} m_{1}(X) m_{2}(Y)$ is called degree of conflict between the sources $\mathcal{B}_1$ and $\mathcal{B}_2$. When $k_{12}=1$,  the orthogonal sum $m (.)$ does not exist and the bodies of evidences $\mathcal{B}_1$ and $\mathcal{B}_2$ are said to be in {\it{full contradiction}}. This rule of combination can be extended easily for the combination of $n>2$ independent sources of evidence. The DST is attractive for the {\it{Data Fusion community}} because it gives a nice mathematical model for ignorance and it includes the Bayesian theory as a special case \cite{Shafer_1976} (p. 4). Although very appealing, the DST presents some weaknesses and limitations because of the Shafer's model itself, the theoretical justification of the Dempster's rule of combination but also because of our confidence to trust the result of Dempster's rule of combination when the conflict becomes important (close to one).\\

The foundations of the DSmT is to abandon the previous Shafer's model (i.e. the exclusivity constraint between elementary hypotheses $\theta_i$ of $\Theta$) just because for some fusion problems it is impossible to define/characterize the problem in terms of well-defined/precise and exclusive elements).  The free-DSm model, denoted $\mathcal{M}^f(\Theta)$, on which is based DSmT allows us to deal with  imprecise/vague notions and concepts between elements of the frame of discernment $\Theta$. The DSmT  includes the possibility to deal with evidences arising from different sources of information which don't have 
access to absolute interpretation of the elements $\Theta$ under consideration and can be interpreted as a general 
and direct extension of probability theory and the Dempster-Shafer theory in the 
following sense. Let $\Theta=\{\theta_{1},\theta_{2}\}$ be the simplest frame of 
discernment involving only two elementary hypotheses (with no more 
additional assumptions on 
$\theta_{1}$ and $\theta_{2}$), then 
\begin{itemize}
\item the probability theory deals with basic probability assignments 
$m(.)\in [0,1]$ such that $m(\theta_{1})+m(\theta_{2})=1$
\item the DST deals with bba $m(.)\in [0,1]$ such that 
$m(\theta_{1})+m(\theta_{2})+m(\theta_{1}\cup\theta_{2})=1$
\item the DSmT theory deals with
new bba $m(.)\in [0,1]$ such that 
$$m(\theta_{1})+m(\theta_{2})+m(\theta_{1}\cup\theta_{2})+m(\theta_{1}\cap\theta_{2})=1$$
\end{itemize}
This comparison can be generalized for the $n$-dimensional case ($|\Theta|=n \geq 2$) as well.
From this very simple idea and from any frame $\Theta$, a new space $D^\Theta$ (free Boolean pre-algebra generated by $\Theta$ and operators $\cap$ and $\cup$), called {\it{hyperpower-set}} is defined \cite{Dezert_2004} as follows:
\begin{enumerate}
\item $\emptyset, \theta_1,\ldots, \theta_n \in D^\Theta$
\item $\forall A\in D^\Theta, B\in D^\Theta, (A\cup B)\in D^\Theta, (A\cap B)\in D^\Theta$
\item No other elements belong to $D^\Theta$, except those, obtained by using rules 1 or 2.
\end{enumerate}
The generation of hyperpower-set $D^\Theta$ is related with the famous Dedekind's problem on enumerating the set of isotone Boolean functions \cite{Dezert_2003f}. From a general frame of discernment $\Theta$ and by adopting the free-DSm model, one then defines a map $m_i(.): D^\Theta \rightarrow [0,1]$, associated to a given source of evidence $\mathcal{B}_i$ such that $m_i(\emptyset)=0$ and $\sum_{A\in D^\Theta} m_i(A) = 1$. This approach allows us to model any source which supports paradoxical (or intrinsic conflicting) information. From this very simple free-DSm model $\mathcal{M}^f(\Theta)$, the classical DSm rule of combination $m(.)\triangleq [m_{1}\oplus \ldots \oplus m_{k}](.)$ of $k\geq 2$ intrinsic conflicting and/or uncertain independent sources of information \begin{equation}
m_{\mathcal{M}^f(\Theta)}(A)=
\sum_{\overset{X_1,\ldots,X_k\in D^\Theta}{X_1\cap \ldots\cap X_k=A}} \prod_{i=1}^{k}m_i(X_i)
\label{eq:DSMClassick}
\end{equation}
\noindent
and with $m_{\mathcal{M}^f(\Theta)}(\emptyset)=0$ by definition. This rule, dealing with uncertain and/or paradoxical/conflicting information is commutative and associative and requires no normalization procedure. A justification of DSm rule of combination from a semantic point of view can be found in \cite{Dezert_2004Book}.\\

Very recently, F. Smarandache and J. Dezert have extended the framework of the DSmT and the previous DSm rule of combination for solving a wider class of fusion problems in which neither free-DSm or Shafer's models fully hold \cite{Dezert_2004Book}. This large class of problems corresponds to problems characterized by any DSm hybrid-model. A DSm hybrid model is defined from the free-DSm model $\mathcal{M}^f(\Theta)$ by introducing some integrity constraints on some elements $A\in D^\Theta$, if there are some certain facts in accordance with the exact nature of the model related to the problem under consideration \cite{Dezert_2004Book}. An integrity constraint on $A\in D^\Theta$ consists in forcing $A$ to be empty through the model $\mathcal{M}$, denoted as $A\overset{\mathcal{M}}{\equiv}\emptyset$. There are several possible kinds of integrity constraints introduced in any free-DSm model:
\begin{itemize}
\item {\it{Exclusivity constraints}}: when some conjunctions of elements of $\Theta$
 are truly impossible, for example when  $\theta_i\cap\ldots\cap\theta_k\overset{\mathcal{M}}{\equiv}\emptyset$.
\item {\it{Non-existential constraints}}: when some disjunctions of elements of $\Theta$ are truly impossible, for example when 
$\theta_i\cup\ldots\cup\theta_k\overset{\mathcal{M}}{\equiv}\emptyset$. The vacuous DSm hybrid model $\mathcal{M}_{\emptyset}$, 
defined by constraint according to the total ignorance: $I_t\triangleq\theta_1\cup \theta_2\cup \ldots\cup \theta_n\overset{\mathcal{M}}{\equiv}\emptyset$, is excluded from consideration, because it is meaningless.
\item {\it{Mixture of exclusivity and non-existential constraints}}: like for example $(\theta_i\cap\theta_j)\cup \theta_k$ or any other hybrid proposition/element of $D^\Theta$ involving both $\cap$ and $\cup$ operators such that at least one element $\theta_k$ is subset of the constrained proposition.
\end{itemize}
The introduction of a given integrity constraint $A\overset{\mathcal{M}}{\equiv}\emptyset \in D^\Theta$ implies the set of inner constraints $B\overset{\mathcal{M}}{\equiv}\emptyset$ for all $B\subset A$. The introduction of two integrity constraints on $A,B \in D^\Theta$ implies the constraint $(A\cup B)\in D^\Theta \equiv\emptyset$ and this implies the emptiness of all $C\in D^\Theta$ such that $C\subset (A\cup B)$. The Shafer's model $\mathcal{M}^{0}(\Theta)$, can be considered as the most constrained DSm hybrid model including all possible exclusivity constraints {\it{without non-existential constraint}}, since all elements in the frame are forced to be mutually exclusive. The DSm hybrid rule of combination, associated to a given DSm hybrid model $\mathcal{M}\neq\mathcal{M}_{\emptyset}$ , for $k\geq 2$ independent sources of information is defined for all $A\in D^\Theta$ as \cite{Dezert_2004Book}:
\begin{equation}
m_{\mathcal{M}(\Theta)}(A)\triangleq 
\phi(A)\Bigl[ S_1(A) + S_2(A) + S_3(A)\Bigr]
 \label{eq:DSmHkBis}
\end{equation}
\noindent
where $\phi(A)$ is the characteristic emptiness function of a set $A$, i.e. $\phi(A)= 1$ if  $A\notin \boldsymbol{\emptyset}$ and $\phi(A)= 0$ otherwise, where $\boldsymbol{\emptyset}\triangleq\{\boldsymbol{\emptyset}_{\mathcal{M}},\emptyset\}$. $\boldsymbol{\emptyset}_{\mathcal{M}}$ is the set  of all elements of $D^\Theta$ which have been forced to be empty through the constraints of the model $\mathcal{M}$ and $\emptyset$ is the classical/universal empty set. $S_1(A)\equiv m_{\mathcal{M}^f(\theta)}(A)$, $S_2(A)$, $S_3(A)$ are defined by \cite{Dezert_2004}
\begin{equation}
S_1(A)\triangleq \sum_{\overset{X_1,X_2,\ldots,X_k\in D^\Theta}{(X_1\cap X_2\cap\ldots\cap X_k)=A}} \prod_{i=1}^{k} m_i(X_i)
\end{equation}
\begin{equation}
S_2(A)\triangleq \sum_{\overset{X_1,X_2,\ldots,X_k\in\boldsymbol{\emptyset}}{ [\mathcal{U}=A]\vee [\mathcal{U}\in\boldsymbol{\emptyset}) \wedge (A=I_t)]}} \prod_{i=1}^{k} m_i(X_i)\end{equation}
\begin{equation}
S_3(A)\triangleq\sum_{\overset{X_1,X_2,\ldots,X_k\in D^\Theta}{\overset{X_1\cup X_2\cup\ldots\cup X_k=A}{\overset{X_1\cap X_2\cap \ldots\cap X_k\in\boldsymbol{\emptyset}}{}}}}  \prod_{i=1}^{k} m_i(X_i)
\end{equation}
with $\mathcal{U}\triangleq u(X_1)\cup u(X_2)\cup \ldots \cup u(X_k)$ where $u(X)$ is the union of all singletons $\theta_i$ that compose $X$ and $I_t \triangleq \theta_1\cup \theta_2\cup\ldots\cup \theta_n$ is the total ignorance. $S_1(A)$ corresponds to the classic DSm rule of combination based on the free-DSm model; $S_2(A)$ represents the mass of all relatively and absolutely empty sets which is transferred to the total or relative ignorances; $S_3(A)$ transfers the sum of relatively empty sets to the non-empty sets.\\ 

{\bf{Theorem}} \cite{Dezert_2004Book}: Let $\Theta = \{\theta_1,\theta_2,\ldots,\theta_n\}$ and it's power-set $2^\Theta=\{\emptyset, \theta_1,\theta_2,\ldots,\theta_n, \text{all possible unions of $\theta_i$'s}\}$, (no intersection), and $k\geq 2$ independent sources of information. After using the DSm classical rule (based on the free-DSm model), one gets belief assignments on many intersections. Let's suppose the Shafer's model $\mathcal{M}^0$, (i.e. all intersections are empty).  Then using next the DSm hybrid rule of combining, one gets:  $m_{\mathcal{M}^0}(X_1\cap\ldots\cap X_k)=0$ and its mass $m(X_1\cap\ldots\cap X_k)$ obtained by the classical DSm rule is transferred to the ignorance element $u(X_1)\cup\ldots\cup u(X_k)$ according to $S_2$ and $S_3$.

\section{First infinite class of counter examples}

The first infinite class of counter examples for the Dempster's rule of combination consists trivially in all cases for which the Dempster's rule becomes mathematically not defined, i.e. one has $0/0$, because of full conflicting sources. The first sub-class presented in subsection \ref{CounterExamples:Sec1:SubSec1.1} corresponds to  Bayesian belief functions. The subsection  \ref{CounterExamples:Sec1:SubSec1.2} will present counter-examples for more general conflicting sources of evidence.

\subsection{Counter-examples for Bayesian sources}
\label{CounterExamples:Sec1:SubSec1.1}

The following examples are devoted only to Bayesian sources, i.e. sources for which the focal elements of belief functions coincide only with some singletons $\theta_i$ of $\Theta$.

\subsubsection{Example with $\Theta=\{\theta_1,\theta_2\}$}

Let's consider the frame of discernment $\Theta=\{\theta_1,\theta_2\}$, two independent experts, and the basic belief masses:
$$m_1(\theta_1) = 1 \qquad m_1(\theta_2) = 0$$
$$m_2(\theta_1) = 0 \qquad m_2(\theta_2) = 1$$

\noindent We represent these belief assignments by the mass matrix 
\begin{equation*}
\mathbf{M}=
\begin{bmatrix}
1 & 0 \\
0 & 1
\end{bmatrix}
\end{equation*}
\begin{itemize}
\item The Dempster's rule can not be applied because one formally gets $m(\theta_1) = 0/0$ and $m(\theta_2) = 0/0$ as well, i.e. undefined.
\item The DSm rule works here because one obtains
$m(\theta_1) = m(\theta_2) = 0$ and $m(\theta_1\cap \theta_2) = 1$ (the total paradox, which really is! if one accepts the free-DSm model). If one adopts the Shafer's model and apply the DSm hybrid rule, then one gets $m_h(\theta_1\cup\theta_2)=1$ which makes sense in this case. The index $h$ denotes here the mass obtained with DSm hybrid rule to avoid confusion with result obtained with DSm classic rule. 
\end{itemize}

\subsubsection{Example with $\Theta=\{\theta_1,\theta_2,\theta_3,\theta_4\}$}

Let's consider the frame of discernment $\Theta=\{\theta_1,\theta_2,\theta_3,\theta_4\}$, two independent experts, and the mass matrix 
\begin{equation*}
\begin{bmatrix}
0.6 & 0  & 0.4 & 0\\
0 & 0.2 & 0 & 0.8
\end{bmatrix}
\end{equation*}
\begin{itemize}
\item Again, the Dempster's rule can not be applied because: $\forall 1\leq j \leq 4$, one gets $m(\theta_j) = 0/0$ (undefined!).
\item But the DSm rule works because one obtains: $m(\theta_1) = m(\theta_2) = m(\theta_3) = m(\theta_4) = 0$, 
and $m(\theta_1 \cap \theta_2) = 0.12$, $m(\theta_1 \cap \theta_4) = 0.48$,
$m(\theta_2 \cap \theta_3) = 0.08$, $m(\theta_3 \cap \theta_4) = 0.32$ (partial paradoxes/conflicts).
\item Suppose now one finds out that all intersections are empty (Shafer's model), then one applies the DSm hybrid rule and one gets (index $h$ stands here for {\it{hybrid}} rule): $m_h(\theta_1\cup \theta_2)=0.12$, $m_h(\theta_1\cup\theta_4)=0.48$, $m_h(\theta_2\cup\theta_3)=0.08$ and $m_h(\theta_3\cup\theta_4)=0.32$.
\end{itemize}

\subsubsection{Another example with $\Theta=\{\theta_1,\theta_2,\theta_3,\theta_4\}$}

LetÕs consider the frame of discernment $\Theta=\{\theta_1,\theta_2,\theta_3,\theta_4\}$, three independent experts, and the mass matrix 
\begin{equation*}
\begin{bmatrix}
0.6 & 0  & 0.4 & 0\\
0 & 0.2 & 0 & 0.8\\
0 & 0.3 & 0 & 0.7
\end{bmatrix}
\end{equation*}
\begin{itemize}
\item Again, the Dempster's rule can not be applied because: $\forall 1\leq j \leq 4$, one gets $m(\theta_j) = 0/0$ (undefined!).
\item But the DSm rule works because one obtains: $m(\theta_1) = m(\theta_2) = m(\theta_3) = m(\theta_4) = 0$, 
and 
$$m(\theta_1\cap\theta_2)=0.6\cdot 0.2 \cdot 0.3 = 0.036$$
$$m(\theta_1 \cap \theta_4) = 0.6\cdot 0.8\cdot 0.7 = 0.336$$
$$m(\theta_2 \cap \theta_3) = 0.4\cdot 0.2\cdot 0.3 = 0.024$$
$$m(\theta_3 \cap \theta_4) = 0.4\cdot 0.8\cdot 0.7 = 0.224$$
$$m(\theta_1\cap\theta_2\cap \theta_4) = 0.6\cdot 0.2\cdot 0.7+0.6\cdot 0.3\cdot 0.8 = 0.228$$
$$m(\theta_2\cap\theta_3\cap \theta_4) = 0.2\cdot 0.4\cdot 0.7+0.3\cdot0.4\cdot 0.8 = 0.152$$
\noindent (partial paradoxes/conflicts) and the others equal zero. If we add all these masses, we get the sum equals to 1.
\item
Suppose now one finds out that all intersections are empty (Shafer's model), then one applies the DSm hybrid rule and one gets: $m_h(\theta_1\cup\theta_2)=0.036$, $m_h(\theta_1\cup\theta_4)=0.336$, $m_h(\theta_2\cup\theta_3)=0.024$, $m_h(\theta_3\cup\theta_4)=0.224$, $m_h(\theta_1\cup\theta_2\cup\theta_4)=0.228$, $m_h(\theta_2\cup\theta_3\cup\theta_4)=0.152$.
\end{itemize}

\subsubsection{More general}

LetÕs consider the frame of discernment $\Theta= \{\theta_1, \theta_2, \ldots, \theta_n\}$, with $n \geq 2$, and $k$ experts, for $k \geq 2$.  Let $\mathbf{M} = [a_{ij}]$, $1\leq i \leq k$, $1\leq j \leq n$, be the mass matrix with $k$ rows and $n$ columns.
If each column of the mass matrix contains at least a zero, then the Dempster's rule can not be applied because one obtains for all $1\leq j \leq n$, $m(\theta_j) = 0/0$ which is undefined!  The degree of conflict is 1.
However, one can use the classical DSm rule and one obtains:
for all $1\leq j \leq n$, $m(\theta_j) = 0$, 
and also partial paradoxes/conflicts: $\forall 1\leq v_s \leq n$, $1 \leq s \leq w$, and $2 \leq w \leq k$, 
$m(\theta_{v_1} \cap \theta_{v_2} \cap \ldots \cap  \theta_{v_w} ) = \sum (a_{1t_1}) \cdot (a_{2t_2})
\cdot \ldots \cdot (a_{kt_k})$, where the set 
$T =\{t_1, t_2, \ldots, t_k\}$ is equal to the set 
$V = \{v_1, v_2, \ldots, v_w\}$ but the order may be different and the elements in the set $T$ could be repeated;  we mean from set $V$ one obtains set $T$ if one repeats some elements of $V$;
therefore: summation $\sum$ is done upon all possible combinations of elements from columns $v_1$, $v_2$, \ldots, $v_w$ such that at least one element one takes from each of these columns $v_1$, $v_2$, \ldots, $v_w$ and also such that from each row one takes one element only;  
the product $(a_{1t_1}) \cdot (a_{2t_2})\cdot \ldots \cdot (a_{kt_k})$ contains one element only from each row $1$, $2$, \ldots, $k$ respectively, and one or more elements from each of the columns $v_1$, $v_2$, \ldots, $v_w$ respectively.

\subsection{Counter-examples for more general sources}
\label{CounterExamples:Sec1:SubSec1.2}

We present in this section two numerical examples involving general (i.e. non Bayesian) sources where the Dempster's rule cannot be applied.

\subsubsection{Example with $\Theta=\{\theta_1,\theta_2,\theta_3,\theta_4\}$}

Let's consider $\Theta=\{\theta_1,\theta_2,\theta_3,\theta_4\}$, two independent experts, and the mass matrix:
\begin{center}
\begin{tabular}{|l|c|c|c|c|c|}
\hline
                   &  $\theta_1$ & $\theta_2$ & $\theta_3$ & $\theta_4$ & $\theta_1\cup\theta_2$\\
\hline
$m_1(.)$  & $0.4$ &  $0.5$       & $0$ &  $0$ & $0.1$\\
$m_2(.)$  & $0$       &  $0$ & $0.3$ &  $0.7$ & $0$\\
\hline
\end{tabular}
\end{center}

The Dempster's rule cannot apply here because one gets $0/0$ for all $m(\theta_i)$, $1\leq i \leq 4$, but the DSm rules (classical or hybrid) work.\\

\noindent Using the DSm classical rule: $m(\theta_1\cap\theta_3) = 0.12$, $m(\theta_1\cap\theta_4)=0.28$, $m(\theta_2\cap\theta_3)=0.15$, $m(\theta_2\cap\theta_4)=0.35$, $m(\theta_3\cap(\theta_1\cup\theta_2))=0.03$, $m(\theta_4\cap(\theta_1\cup\theta_2))=0.07$.\\

\noindent
Suppose now one finds out that one has a Shafer's model; then one uses the 
DSm hybrid rule (denoted here with index $h$): $m_h(\theta_1\cup\theta_3) = 0.12$, $m_h(\theta_1\cup\theta_4)=0.28$, $m_h(\theta_2\cup\theta_3)=0.15$, 
$m_h(\theta_2\cup\theta_4)=0.35$, $m_h(\theta_3\cup\theta_1\cup\theta_2)=0.03$, $m_h(\theta_4\cup\theta_1\cup\theta_2)=0.07$.

\subsubsection{Another example with $\Theta=\{\theta_1,\theta_2,\theta_3,\theta_4\}$}
Let's consider $\Theta=\{\theta_1,\theta_2,\theta_3,\theta_4\}$, three independent experts, and the mass matrix:
\begin{center}
\begin{tabular}{|l|c|c|c|c|c|c|}
\hline
                  &  $\theta_1$ & $\theta_2$ & $\theta_3$ & $\theta_4$ & $\theta_1\cup\theta_2$ & $\theta_3\cup\theta_4$\\
\hline
$m_1(.)$  & $0.4$ &  $0.5$       & $0$ &  $0$ & $0.1$ & $0$\\
$m_2(.)$  & $0$       &  $0$ & $0.3$ &  $0.6$ & $0$& $0.1$\\
$m_3(.)$  & $0.8$       &  $0$ & $0$ &  $0$ & $0.2$& $0$\\
\hline
\end{tabular}
\end{center}

The Dempster's rule cannot apply here because one gets $0/0$ for all $m(\theta_i)$, $1\leq i \leq 4$, but the DSm rules (classical or hybrid) work.\\

\noindent 
Using the DSm classical rule, one gets:
\begin{center}
\begin{tabular}{ll}
$m(\theta_1)=m(\theta_2)=m(\theta_3)=m(\theta_4)=0$ &  $m(\theta_1\cup\theta_2)=m(\theta_3\cup\theta_4)=0$\\
$m(\theta_1\cap\theta_3)=0.096$ & $m(\theta_1\cap\theta_3\cap(\theta_1\cup\theta_2))=m(\theta_1\cap\theta_3)=0.024$\\
$m(\theta_1\cap\theta_4)=0.192$ &$m(\theta_1\cap\theta_4\cap(\theta_1\cup\theta_2))=m(\theta_1\cap\theta_4)=0.048$\\
$m(\theta_1\cap(\theta_3\cup\theta_4))=0.032$ & $m(\theta_1\cap(\theta_3\cup\theta_4)\cap(\theta_1\cup\theta_2))=m(\theta_1\cap(\theta_3\cup\theta_4))=0.008$\\
$m(\theta_2\cap\theta_3\cap\theta_1)=0.120$ & $m(\theta_2\cap\theta_3\cap(\theta_1\cup\theta_2))=m(\theta_2\cap\theta_3)=0.030$\\
$m(\theta_2\cap\theta_4\cap\theta_1)=0.240$ & $m(\theta_2\cap\theta_4\cap(\theta_1\cup\theta_2))=m(\theta_2\cap\theta_4)=0.060$\\
$m(\theta_2\cap(\theta_3\cup\theta_4)\cap\theta_1)=m((\theta_1\cap\theta_2)\cap(\theta_3\cup\theta_4))=0.040$ & $m(\theta_2\cap(\theta_3\cup\theta_4)\cap(\theta_1\cup\theta_2))=m(\theta_2\cap(\theta_3\cup\theta_4))=0.010$\\
$m((\theta_1\cup\theta_2)\cap\theta_3\cap\theta_1)=m(\theta_1\cap\theta_3)=0.024$ & $m((\theta_1\cup\theta_2)\cap\theta_3)=0.006$\\
$m((\theta_1\cup\theta_2)\cap\theta_4\cap\theta_1)=m(\theta_1\cap\theta_4)=0.048$ & $m((\theta_1\cup\theta_2)\cap\theta_4)=0.012$\\
$m((\theta_1\cup\theta_2)\cap(\theta_3\cup\theta_4)\cap\theta_1)=m(\theta_1\cap(\theta_3\cup\theta_4))=0.008$ &
$m((\theta_1\cup\theta_2)\cap(\theta_3\cup\theta_4))=0.002$\\
\end{tabular}
\end{center}

\noindent
After cummulating, one finally gets with DSm classic rule:
\begin{center}
\begin{tabular}{ll}
$m(\theta_1\cap\theta_3)=0.096+0.024+0.024=0.144$& $m(\theta_1\cap\theta_4)=0.192+0.048+0.048=0.288$\\
$m(\theta_2\cap\theta_3)=0.030$ & $m(\theta_2\cap\theta_4)=0.060$\\
$m(\theta_1\cap\theta_2\cap\theta_3)=0.120$ & $m(\theta_1\cap\theta_2\cap\theta_4)=0.240$\\
$m((\theta_1\cup\theta_2)\cap\theta_3)=0.006$ & $m((\theta_1\cup\theta_2)\cap\theta_4)=0.012$\\
$m(\theta_1\cap(\theta_3\cup\theta_4))=0.032+0.008+0.008=0.048$ &
$m(\theta_1\cap\theta_2\cap(\theta_3\cup\theta_4))=0.040$\\
$m(\theta_2\cap(\theta_3\cup\theta_4))=0.010$ &
$m((\theta_1\cup\theta_2)\cap(\theta_3\cup\theta_4))=0.002$
\end{tabular}
\end{center}

\noindent
Suppose now, one finds out that all intersections are empty. Using the DSm hybrid rule one gets:

\begin{center}
\begin{tabular}{ll}
$m_h(\theta_1\cup\theta_3)=0.144$ & $m_h(\theta_1\cup\theta_4)=0.288$\\
$m_h(\theta_2\cup\theta_3)=0.030$ & $m_h(\theta_2\cup\theta_4)=0.060$\\
$m_h(\theta_1\cup\theta_2\cup\theta_3)=0.120+0.006=0.126$ & $m_h(\theta_1\cup\theta_2\cup\theta_4)=0.240+0.012=0.252$\\
$m_h(\theta_1\cup\theta_3\cup\theta_4)=0.048$ & $m_h(\theta_2\cup\theta_3\cup\theta_4)=0.010$\\
$m_h(\theta_1\cup\theta_2\cup\theta_3\cup\theta_4)=0.040+0.002=0.042$ &
\end{tabular}
\end{center}

\subsubsection{More general}

LetÕs consider the frame of discernment $\Theta= \{\theta_1, \theta_2, \ldots, \theta_n\}$, with $n \geq 2$, and $k$ experts, for $k \geq 2$, and the mass matrix $\mathbf{M}$ with $k$ rows and $n+u$ columns, where $u \geq 1$, corresponding to $\theta_1$, $\theta_2$, \ldots, $\theta_n$,  and $u$ uncertainties $\theta_{i_1} \cup \ldots \cup \theta_{i_s}$, \ldots, $\theta_{j_1} \cup \ldots \cup \theta_{j_t}$ respectively.\\

\noindent If the following conditions occur:
\begin{itemize}
\item
each column contains at least one zero;
\item
all uncertainties are different from the total ignorance $\theta_{1} \cup \ldots \cup \theta_{n}$ (i.e., they are partial ignorances);
\item
the partial uncertainties are disjoint two by two;
\item
for each non-null uncertainty column $c_j$, $n+1 \leq j \leq n+u$, of the form say $\theta_{p_1} \cup \ldots \cup \theta_{p_w}$, there exists a row such that all its elements on columns $p_1$, \ldots, $p_w$, and $c_j$ are zero.  
\end{itemize}
\noindent
then the Dempster's rule of combination cannot apply for such infinite class of fusion problems because one gets $0/0$ for all $m(\theta_i)$, $1\leq i \leq n$. The DSm rules (classical or hybrid) work for such infinite class of examples.

\section{Second infinite class of counter examples}

This second class of counter-examples generalizes the famous Zadeh's example given in \cite{Zadeh_1979,Zadeh_1986}.

\subsection{Zadeh's example}

Two doctors examine a patient and agree that it suffers from either meningitis 
(M), concussion (C) or brain tumor (T). Thus $\Theta=\{M,C,T\}$. Assume that the doctors 
agree in their low expectation of a tumor, but disagree in likely 
cause and provide the following diagnosis
 $$m_{1}(M)=0.99\qquad m_1(T)=0.01 \qquad\text{and}\qquad m_{2}(C)=0.99\qquad m_2(T)=0.01$$
If we combine the two basic belief functions using Dempster's rule of combination, one 
gets the unexpected final conclusion 
$$m(T)=\frac{0.0001}{1-0.0099-0.0099-0.9801}=1$$
\noindent
 which means that the patient 
suffers with certainty from brain tumor !!!. This unexpected 
result arises from the fact that the two bodies of evidence (doctors) agree 
that patient does not suffer from tumor but are in almost full 
contradiction for the other causes of the disease. This very simple 
but interesting example shows the limitations of practical use of the DST for automated 
reasoning.\\

This example has been examined in literature by several authors to explain the anomaly 
of the result of  the Dempster's rule of combination in such case. Due to the high degree of conflict arising in such extreme case willingly pointed out by Zadeh to show the weakness of this rule, it is often argued that in such case the result of the Dempster's rule must not be taken directly without checking the level of the conflict between sources of evidence. This is trivially true but there is no theoretical way to decide {\it{beforehand}} if one can trust or not the result of such rule of combination. This is one of its major drawback. The issue consists generally in choosing rather somewhat arbitrarily or heuristically some threshold value on the degree of conflict between sources to accept or reject the result of the fusion. Such approach can't be solidly justified from theoretical analysis. Assuming such threshold is set to a given value, say 0.70 for instance, is it acceptable to reject the fusion result if the conflict appears to be 0.7001 and accept it when the conflict becomes 0.6999? What to do when the decision about the fusion result is rejected and one has no assessment on the reliability of the sources or when the sources have the same reliability/confidence but an important decision has to be taken anyway? There is no theoretical solid justification which can reasonably support such kind of approaches commonly used in practice.\\

The two major explanations of this problem found in literature are mainly based either on the fact that problem arises from the closed-world assumption of the Shafer's model $\Theta$ and it is suggested to work rather with an open-world model and/or the fact that sources of evidence are not reliable. These explanations although being admissible are not necessarily the only correct explanations.  Note that the open-world assumption can always be easily relaxed advantageously by introducing a new hypothesis, say $\theta_0$ in the initial frame $\Theta=\{\theta_1,\ldots, \theta_n\}$ in order to close it. $\theta_0$ will then represent all possible alternatives (although remaining unknown) of initial hypotheses $\theta_1$,\ldots $\theta_n$. This idea has been already proposed by Yager in \cite{Yager_1983} through his hedging solution. Upon our analysis, it is not necessary to adopt/follow the open-world model neither to admit the assumption about the reliability of the sources to find a justification in this counter-intuitive result. Actually, both sources can have the same reliability and the Shafer's model can be accepted for the combination of the two reports by using another rule of combination. This is exactly the purpose of the DSm hybrid rule of combination. Of course when one has some prior information on the reliability of sources, one has to take them into account properly by some discounting methods. The discounting techniques can also apply in the DSmT framework and there is no incompatibility to mix both discounting techniques with DSm rules of combinations when necessary (when there is strong reason to justify doing it, i.e. when one has prior reliable information on reliability of the sources). The discounting techniques must never been used as an artificial ad-hoc mechanism to update Dempster's result once problem has arisen. We strongly disagree with the idea that all problems with the Dempster's rule can be solved {\it{beforehand}} by discounting techniques. This can help obviously to improve the assessment of belief function to be combined when used properly and fairly, but this does not fundamentally solve the inherent problem of the Dempster's rule itself when conflict remains high.\\

The problem comes from the fact that both sources provide essentially their belief with respect only to their own limited knowledge and experience. It is also possible in some cases, that sources of information even don't have the same interpretation of concepts included in the frame of the problem (such kind of situation frequently appears by example in debates on TV, on radio or in most of meetings where important decision/approval have to be drawn and when the sources don't share the same opinion. This is what happens daily in real life and one has to deal with such conflicting situations anyway). In this Zadeh's controversy example, it is possible that the first doctor is expert mainly in meningitis and in brain tumor while the second doctor is expert mainly in cerebral contusion and in brain tumor. Because of their limited knowledges and experiences, both doctors can also have also the same reliability. If they have been asked to give their reports only on $\Theta=\{M,C,T\}$ (but not on an extended frame), their reports have to be taken with same weight and the combination have to be done anyway when one has no solid reason to reject one report with respect to the other one; the result of the Demsper's rule still remains very questionable. No rational brain surgeon would  take the decision for a brain intervention (i.e. a risky tumor ablation) based on the Dempster's rule result, neither the family of the patient. Therefore upon our analysis, the two previous explanations given in literature (although being possible and admissible in some cases) are not necessary and sufficient to explain the source of the anomaly. Several alternatives to the Dempster's rule to circumvent this anomaly have been proposed in literature mainly through the works of R. Yager  \cite{Yager_1983}, D. Dubois and H. Prade  \cite{Dubois_1986} already reported in \cite{Dezert_2004Book} or by Daniel in \cite{Daniel_2000}. The DSmT offers just a new issue for solving also such controversy example as it will be showed. In summary, some extreme caution on the degree of conflict of the sources must always be taken before taking a final decision 
based on the Dempster's rule of combination, specially when vital wagers are involved.\\

If we now adopt the free-DSm model, i.e. we replace the initial Shafer's model by accepting the possibility of non null intersection between hypotheses $M$, $C$ and $T$ and by working directly on hyper-power set $D^\Theta$ then one gets directly and easily the following result with the classical DSm rule of combination:
$$m(M\cap C)=0.9801 \qquad m(M\cap T)=0.0099 \qquad m(C\cap T)=0.0099 \qquad m(T)=0.0001$$
\noindent which makes sense when working with such new model. Obviously same result can be obtained (the proof is left here to the reader) when working with the Dempster's rule based on the following refined frame $\Theta_{ref}$ defined with basic belief functions on power-set $2^{\Theta_{ref}}$:
\begin{equation*}
\begin{split}
\Theta_{ref}=\{\theta_1=M\cap C\cap T,\theta_2=M\cap C\cap \bar{T}, \theta_3=M\cap \bar{C}\cap T, \theta_4=\bar{M}\cap C\cap T,\\
\theta_5=M\cap \bar{C}\cap \bar{T},
\theta_6=\bar{M}\cap C \cap \bar{T},\theta_7=\bar{M}\cap \bar{C} \cap T\}
\end{split}
\end{equation*}
\noindent where $\bar{T}$,$\bar{C}$ and $\bar{M}$ denote respectively the complement of $T$, $C$ and $M$.\\

The equality of both results (i.e. by the classical DSm rule based on the free-DSm model and by the Dempster's rule based on the refined frame) is just normal since the normalization factor $1-k$ of the Dempster's rule in this case reduces to 1 because of the new choice of the new model. Based on this remark, one could then try to argue that DSmT (together with its DSm classical rule for free-DSm model) is superfluous. Such claim is obviously wrong for the two following reasons: it is unecessary to work with a bigger space (keeping in mind that $|D^\Theta| < |2^{\Theta_{ref}}|$) to get the result (the DSm rule offers just a direct and more convenient issue to get the result), but also because in some fusion problems involving vague/continuous concepts, the refinement is just impossible to obtain and we are unfortunately forced to deal with ambiguous concepts/hypotheses (see \cite{Goodman_1997} for details and justification).\\

If one has no doubt on the reliability of both Doctors (or no way to assess it) and if one is absolutely sure that the true origin of the suffering of the patient lies only in the frame $\Theta=\{M,C,T\}$ and we consider these origins as truly exclusive, then one has to work with the initital frame of discernment $\Theta$ satisfying the Shafer's model. As previously showed, the Dempster's rule fails to provide a reasonable and acceptable conclusion in such high conflicting case. However, this case can be easily handled by the DSm hybrid rule of combination.  DSm hybrid rule applies now because the Shafer's model is nothing but a  particular hybrid model including all exclusivity constraints between hypotheses of the frame $\Theta$ \cite{Dezert_2004Book}. One then gets with the DSm hybrid rule for this simple case (more general and complex examples are given in \cite{Dezert_2004Book}), after the proper mass transfer of all sources of the conflicts:
$$m(M\cup C)=0.9801 \qquad m(M\cup T)=0.0099 \qquad m(C\cup T)=0.0099 \qquad m(T)=0.0001$$
This result is not surprising and makes perfectly sense with common intuition actually since it provides a coherent and reasonable solution to the problem. It shows clearly that a brain intervention for ablation of an hypothetical tumor is not recommended, but preferentially a better examination of the patient focused on Meningitis or Contusion as possible source of the suffering. The consequence of the results of the Dempster's rule and the DSm hybrid rule is therefore totally different. 

\subsection{Generalization with $\Theta=\{\theta_1,\theta_2,\theta_3\}$}

Let's consider $0 < \epsilon_1,\epsilon_2 < 1$ be two very tiny positive numbers (close to zero), the frame of discernment be $\Theta=\{\theta_1,\theta_2,\theta_3\}$, have two experts (independent sources of evidence $s_1$ and $s_2$) giving the belief masses
$$m_1(\theta_1)=1-\epsilon_1 \quad m_1(\theta_2)=0 \quad m_1(\theta_3)=\epsilon_1$$
$$m_2(\theta_1)=0 \quad m_2(\theta_2)=1-\epsilon_2 \quad m_2(\theta_3)=\epsilon_2$$
\noindent From now on, we prefer to use matrices to describe the masses, i.e.
$$\begin{bmatrix}
1-\epsilon_1 & 0 &\epsilon_1\\
0 & 1-\epsilon_2 & \epsilon_2
\end{bmatrix}
$$
\begin{itemize}
\item Using the Dempster's rule of combination, one gets
$$m(\theta_3)=\frac{(\epsilon_1\epsilon_2)}{(1-\epsilon_1)\cdot 0 + 0\cdot (1-\epsilon_2) + \epsilon_1\epsilon_2}=1$$
\noindent which is absurd (or at least counter-intuitive). Note that whatever positive values for $\epsilon_1$, $\epsilon_2$ are, the Dempster's rule of combination provides always the same result (one) which is abnormal. The only acceptable and correct result obtained by the Dempster's rule is really obtained only in the trivial case when $\epsilon_1=\epsilon_2=1$, i.e. when both sources agree in $\theta_3$ with certainty which is obvious.
\item Using the DSm rule of combination based on free-DSm model, one gets 
$m(\theta_3)=\epsilon_1\epsilon_2$, $m(\theta_1\cap \theta_2)=(1-\epsilon_1)(1-\epsilon_2)$, $m(\theta_1\cap \theta_3)=(1-\epsilon_1)\epsilon_2$, $m(\theta_2\cap \theta_3)=(1-\epsilon_2)\epsilon_1$ and the others are zero which appears more reliable/trustable.
\item Go back to the Shafer's model and using the DSm hybrid rule of combination, one gets 
$m(\theta_3)=\epsilon_1\epsilon_2$, $m(\theta_1\cup \theta_2)=(1-\epsilon_1)(1-\epsilon_2)$, $m(\theta_1\cup \theta_3)=(1-\epsilon_1)\epsilon_2$, $m(\theta_2\cup \theta_3)=(1-\epsilon_2)\epsilon_1$ and the others are zero.
\end{itemize}

\noindent
Note that in the special  case when $\epsilon_1=\epsilon_2=1/2$, one has
$$m_1(\theta_1)=1/2 \quad m_1(\theta_2)=0 \quad m_1(\theta_3)=1/2 \qquad\text{and}\qquad m_2(\theta_1)=0 \quad m_2(\theta_2)=1/2 \quad m_2(\theta_3)=1/2$$
The Dempster's rule of combinations still yields $m(\theta_3)=1$ while the DSm hybrid rule based on the same Shafer's model yields now
$m(\theta_3)=1/4$, $m(\theta_1\cup \theta_2)=1/4$, $m(\theta_1\cup \theta_3)=1/4$, $m(\theta_2\cup \theta_3)=1/4$ which is normal.

\subsection{Generalization with $\Theta=\{\theta_1,\theta_2,\theta_3,\theta_4\}$}

Let's consider $0< \epsilon_1,\epsilon_2,\epsilon_3 < 1$ be three very tiny positive numbers, the frame of discernment be $\Theta=\{\theta_1,\theta_2,\theta_3,\theta_4\}$, have two experts giving the mass matrix
$$\begin{bmatrix}
1-\epsilon_1-\epsilon_2 & 0 &\epsilon_1& \epsilon_2\\
0 & 1-\epsilon_3& 0 & \epsilon_3
\end{bmatrix}
$$
Again using the Dempster's rule of combination, one gets $m(\theta_4)=1$ which is absurd while using the DSm rule of combination based on free-DSm model, one gets  $m(\theta_4)=\epsilon_2\epsilon_3$ which is reliable. Using the DSm classical rule:
$m(\theta_1\cap\theta_2)=(1-\epsilon_1-\epsilon_2)(1-\epsilon_3)$, $m(\theta_1\cap\theta_4)=(1-\epsilon_1-\epsilon_3)\epsilon_3$, 
$m(\theta_3\cap\theta_2)=\epsilon_1(1-\epsilon_3), m(\theta_3\cap\theta_4)=\epsilon_1\epsilon_3$,
$m(\theta_4)=\epsilon_2\epsilon_3$. Suppose one finds out that all intersections are empty, then one applies the DSm hybrid rule:
$m_h(\theta_1\cup\theta_2)=(1-\epsilon_1-\epsilon_2)(1-\epsilon_3)$, $m_h(\theta_1\cup\theta_4)=(1-\epsilon_1-\epsilon_3)\epsilon_3$, 
$m_h(\theta_3\cup\theta_2)=\epsilon_1(1-\epsilon_3)$, $m_h(\theta_3\cup\theta_4)=\epsilon_1\epsilon_3$,
$m_h(\theta_4)=\epsilon_2\epsilon_3$.

\subsection{More general}

Let's consider $0< \epsilon_1,\ldots,\epsilon_n < 1$ be very tiny positive numbers, the frame of discernment be $\Theta=\{\theta_1,\ldots,\theta_n,\theta_{n+1}\}$, have two experts giving the mass matrix
\begin{equation*}
\begin{bmatrix}
1- S_1^p & 0 &\epsilon_1& 0 & \epsilon_2 & \ldots & 0 & \epsilon_p\\
0 & 1- S_{p+1}^n  & 0 &\epsilon_{p+1} & 0 &  \ldots & \epsilon_{n-1} & \epsilon_n
\end{bmatrix}
\end{equation*}
\noindent where $1\leq p \leq n$ and $S_1^p\triangleq\sum_{i=1}^p \epsilon_i$ and $S_{p+1}^n\triangleq\sum_{i=p+1}^n \epsilon_i$. Again using the Dempster's rule of combination, one gets $m(\theta_{n+1})=1$ which is absurd while using the DSm rule of combination based on free-DSm model, one gets  $m(\theta_{n+1})=\epsilon_p\epsilon_n$ which is reliable. This example is similar to the previous one, but generalized.

\subsection{Even more general}

Let's consider $0< \epsilon_1,\ldots,\epsilon_n < 1$ be very tiny positive numbers (close to zero), the frame of discernment be $\Theta=\{\theta_1,\ldots,\theta_n,\theta_{n+1}\}$, have $k\geq 2$ experts giving the mass matrix of $k$ rows and $n+1$ columns such that:
\begin{itemize}
\item
one column, say column $j$, is $(\epsilon_{j_1},\epsilon_{j_2},\ldots,\epsilon_{j_k})'$ (transposed vector), where $1\leq j\leq n+1$ where $\{\epsilon_{j_1},\epsilon_{j_2},\ldots,\epsilon_{j_k}\}$ is included in $\{\epsilon_{1},\epsilon_{2},\ldots,\epsilon_{n}\}$;
\item
and each column (except column $j$) contains at least one element equals to zero.
\end{itemize}
\noindent Then the Dempster's rule of combination gives $m(\theta_j)=1$ which is absurd, while the classical DSm rule gives $m(\theta_j)=\epsilon_{j_1}\cdot \epsilon_{j_2}\cdot \ldots\cdot \epsilon_{j_k}\neq 0$ which is reliable.\\

Actually, we need to set restrictions only for $\epsilon_{j_1}$, $\epsilon_{j_2}$, \ldots, and $\epsilon_{j_k}$ to be very tiny positive numbers, not for all $\epsilon_1$, $\epsilon_2$, \ldots, $\epsilon_n$ (the others can be anything in the interval $[0, 1)$ such that the sum of elements on each row be equal 1).

\section{Third infinite class of counter examples}

This third class of counter-examples deals with belief functions committing a non null mass to some uncertainties.

\subsection{Example with $\Theta=\{\theta_1,\theta_2,\theta_3,\theta_4\}$}

Let's consider $\Theta=\{\theta_1,\theta_2,\theta_3,\theta_4\}$, two independent experts, and the mass matrix:
\begin{center}
\begin{tabular}{|l|c|c|c|c|c|}
\hline
                   &  $\theta_1$ & $\theta_2$ & $\theta_3$ & $\theta_4$ & $\theta_3\cup\theta_4$\\
\hline
$m_1(.)$  & $0.99$ &  $0$       & $0$ &  $0$ & $0.01$\\
$m_2(.)$  & $0$       &  $0.98$ & $0$ &  $0$ & $0.02$\\
\hline
\end{tabular}
\end{center}

If one applies the Dempster's rule, one gets 
$$m(\theta_3\cup\theta_4)=\frac{(0.01\cdot 0.02)}{(0+0+0+0+0.01\cdot 0.02)}=1$$
\noindent (total ignorance), which doesn't bring any information to the fusion. This example looks similar to Zadeh's example, but is different because it is referring to uncertainty (not to contradictory) result.
Using the DSm classical rule:
$m(\theta_1\cap\theta_2)=0.9702$, $m(\theta_1\cap(\theta_3\cup\theta_4))=0.0198$, 
$m(\theta_2\cap(\theta_3\cup\theta_4))=0.0098$, $m(\theta_3\cup\theta_4)=0.0002$.
Suppose now one finds out that all intersections are empty (i.e. one adopts the Shafer's model).
Using the DSm hybrid rule one gets:
$m_h(\theta_1\cup\theta_2)=0.9702$, $m_h(\theta_1\cup\theta_3\cup\theta_4)=0.0198$, $m_h(\theta_2\cup\theta_3\cup\theta_4)=0.0098$, 
$m_h(\theta_3\cup\theta_4)=0.0002$.

\subsection{Example with $\Theta=\{\theta_1,\theta_2,\theta_3,\theta_4,\theta_5\}$}

Let's consider $\Theta=\{\theta_1,\theta_2,\theta_3,\theta_4,,\theta_5\}$, three independent experts, and the mass matrix:
\begin{center}
\begin{tabular}{|l|c|c|c|c|c|c|}
\hline
                   &  $\theta_1$ & $\theta_2$ & $\theta_3$ & $\theta_4$ & $\theta_5$ & $\theta_4\cup\theta_5$\\
\hline
$m_1(.)$  & $0.99$ &  $0$       & $0$       &  $0$ &  $0$ & $0.01$\\
$m_2(.)$  & $0$       &  $0.98$ & $0.01$ &  $0$ &  $0$ & $0.01$\\
$m_3(.)$  & $0.01$ &  $0.01$ & $0.97$ &  $0$ &  $0$ & $0.01$\\
\hline
\end{tabular}
\end{center}

\begin{itemize}
\item
If one applies the Dempster's rule, one gets 
$$m(\theta_4\cup\theta_5)=\frac{(0.01\cdot 0.01\cdot 0.01)}{(0+0+0+0+0.01\cdot 0.01\cdot 0.01)}=1$$
\noindent (total ignorance), which doesn't bring any information to the fusion. 
\item
Using the DSm classical rule one gets:
\begin{equation*}
m(\theta_1\cap\theta_2)=0.99\cdot 0.98\cdot 0.01 +  0.99\cdot 0.98\cdot 0.01= 0.019404
\end{equation*}
\begin{equation*}
m(\theta_1\cap\theta_3)=0.99\cdot 0.01\cdot 0.01+ 0.99\cdot 0.01\cdot 0.97= 0.009702
\end{equation*}
\begin{equation*}
m(\theta_1\cap\theta_2\cap\theta_3)=0.99\cdot 0.98\cdot 0.97+ 0.99\cdot 0.01\cdot 0.01= 0.941193
\end{equation*}
\begin{equation*}
m(\theta_1\cap\theta_3\cap(\theta_4\cup\theta_5)) = 0.99\cdot 0.01\cdot 0.01 + 0.99\cdot 0.01\cdot 0.97 + 0.01\cdot 0.01\cdot 0.01= 0.009703
\end{equation*}
\begin{equation*}
m(\theta_1\cap(\theta_4\cup\theta_5))=0.99\cdot 0.01\cdot 0.01 + 
0.99\cdot 0.01\cdot 0.01 + 0.01\cdot 0.01\cdot 0.01=0.000199
\end{equation*}
\begin{equation*}
m((\theta_4\cup\theta_5)\cap\theta_2\cap\theta_1)=0.01\cdot 0.98\cdot 0.01 + 
0.99\cdot 0.01\cdot 0.01 + 0.99\cdot 0.98\cdot 0.01= 0.009899
\end{equation*}
\begin{equation*}
m((\theta_4\cup\theta_5)\cap\theta_2)=0.01\cdot 0.98\cdot 0.01 + 
 0.01\cdot 0.98\cdot 0.01 + 0.01\cdot 0.01\cdot 0.01= 0.000197
\end{equation*}
\begin{equation*}
m((\theta_4\cup\theta_5)\cap\theta_2\cap\theta_3)=0.01\cdot 0.98\cdot 0.97 + 
0.01\cdot 0.01\cdot 0.01=0.009507
\end{equation*}
\begin{equation*}
m((\theta_4\cup\theta_5)\cap\theta_3)=0.01\cdot 0.01\cdot 0.97+ 
0.01\cdot 0.01\cdot 0.01 + 0.01\cdot 0.01\cdot 0.97= 0.000195
\end{equation*}
\begin{equation*}
m(\theta_4\cup\theta_5)=0.01\cdot 0.01\cdot 0.01= 0.000001
\end{equation*}
\noindent The sum of all masses is 1.
\item
Suppose now one finds out that all intersections are empty (Shafer's model), then one uses the DSm hybrid rule and one gets:
\begin{center}
\begin{tabular}{ll}
$m_h(\theta_1\cup\theta_2)= 0.019404$ & $m_h(\theta_1\cup\theta_3)= 0.009702$\\
$m_h(\theta_1\cup\theta_2\cup\theta_3)= 0.941193$ &$m_h(\theta_1\cup\theta_3\cup\theta_4\cup\theta_5) = 0.009703$\\
$m_h(\theta_1\cup\theta_4\cup\theta_5)=0.000199$ & $m_h(\theta_4\cup\theta_5\cup\theta_2\cup\theta_1)= 0.009899$\\
$m_h(\theta_4\cup\theta_5\cup\theta_2)= 0.000197$ &$m_h(\theta_4\cup\theta_5\cup\theta_2\cup\theta_3)=0.009507$\\
$m_h(\theta_4\cup\theta_5\cup\theta_3)= 0.000195$ & $m_h(\theta_4\cup\theta_5)= 0.000001$
\end{tabular}
\end{center}
\noindent The sum of all masses is 1.
\end{itemize}

\subsection{More general}

Let $\Theta=\{\theta_1,\ldots,\theta_n\}$, where $n \geq 2$, $k$ independent experts, $k \geq 2$, and the mass matrix $\mathbf{M}$ of $k$ rows and $n+1$ columns, corresponding to $\theta_1$, $\theta_2$, \ldots, $\theta_n$,  and one uncertainty (different from the total uncertainty $\theta_1\cup\theta_2\cup\ldots\cup\theta_n$) say $\theta_{i_1}\cup\ldots\cup\theta_{i_s}$ respectively.  If the following conditions occur:
\begin{itemize}
\item
each column contains at least one zero, except the last column (of uncertainties) which has only non-null elements, 
$0 < \epsilon_1, \epsilon_2, \ldots, \epsilon_k < 1$, very tiny numbers (close to zero);
\item
the columns corresponding to the elements $\theta_{i_1}$,\ldots, $\theta_{i_s}$ are null (all their elements are equal to zero).
\end{itemize}
\noindent
If one applies the Dempster's rule, one gets $m(\theta_{i_1}\cup\ldots\cup\theta_{i_s}) = 1$ (total ignorance), which doesnÕt bring any information to the fusion.

\subsection{Even more general}

One can extend the previous case even more, considering to $u$ uncertainty columns, $u \geq 1$ as follows.\\

Let $\Theta=\{\theta_1,\ldots,\theta_n\}$, where $n \geq 2$, $k$ independent experts, $k \geq 2$, and the mass matrix $\mathbf{M}$ of $k$ rows and $n+u$ columns, corresponding to $\theta_1$, $\theta_2$, \ldots, $\theta_n$,  and $u$ uncertainty columns (different from the total uncertainty $\theta_1\cup\theta_2\cup\ldots\cup\theta_n$) respectively.  If the following conditions occur:
\begin{itemize}
\item
each column contains at least one zero, except one column among the last $u$ uncertainty ones which has only non-null elements
$0 < \epsilon_1, \epsilon_2, \ldots, \epsilon_k < 1$, very tiny numbers (close to zero);
\item
the columns corresponding to all elements $\theta_{i_1}$,\ldots, $\theta_{i_s}$,\ldots, $\theta_{r_1}$,\ldots, $\theta_{r_s}$ (of course, these elements should not be all $\theta_1$, $\theta_2$,\ldots, $\theta_n$, but only a part of them) that occur in all uncertainties are null (i.e., all their elements are equal to zero).
\end{itemize}
\noindent
If one applies the Dempster's rule, one gets $m(\theta_{i_1}\cup\ldots\cup\theta_{i_s}) = 1$ (total ignorance), which doesn't bring any information to the fusion.

\section{Fourth infinite class of counter examples}

This infinite class of counter-examples concerns the Dempster's rule of conditioning defined as \cite{Shafer_1976} :

$$\forall B\in 2^\Theta, \quad m(B|A)= \frac{\sum_{X,Y\in 2^\Theta, (X\cap Y)=B} m(X)m_A(Y)}{1-\sum_{X,Y\in 2^\Theta, (X\cap Y)=\emptyset} m(X)m_A(Y)}$$

\noindent 
where $m(.)$ is any proper basic belief function defined over $2^\Theta$ and $m_A(.)$ is a particular belief function defined by choosing $m_A(A)=1$ for any $A\in 2^\Theta$ with $A\neq\emptyset$.

\subsection{Example with $\Theta=\{\theta_1,\ldots,\theta_6\}$}

Let's consider $\Theta=\{\theta_1,\ldots,\theta_6\}$, one expert and a certain body of evidence over $\theta_2$, with the mass matrix:
\begin{center}
\begin{tabular}{|l|c|c|c|c|c|}
\hline
                   &  $\theta_1$ & $\theta_2$ & $\theta_3$ & $\theta_4\cup\theta_5$ & $\theta_5\cup\theta_6$\\
\hline
$m_1(.)$  & $0.3$ &  $0$ &  $0.4$ &  $0.2$ &  $0.1$\\
$m_{\theta_2}(.)$  & $0$ &  $1$ &  $0$ &  $0$ &  $0$\\
\hline
\end{tabular}
\end{center}

\begin{itemize}
\item Using the Dempster's rule of conditioning, one gets: $m(. | \theta_2) = 0/0$ for all the masses.
\item Using the DSm classical rule, one gets: 
%
$$m(\theta_1\cap\theta_2 | \theta_2) = 0.3\qquad m(\theta_2\cap\theta_3 | \theta_2) = 0.4\qquad m(\theta_2\cap(\theta_4\cup\theta_5) | \theta_2) = 0.2\qquad m(\theta_2\cap(\theta_5\cup\theta_6) | \theta_2) = 0.1$$
\item If now, one finds out that all intersections are empty (we adopt the Shafer's model), then using the DSm hybrid rule, one gets:
$$m_h(\theta_1\cup\theta_2|\theta_2)=0.3\qquad m_h(\theta_2\cup\theta_3|\theta_2)=0.4\qquad m_h(\theta_2\cup\theta_4\cup\theta_5|\theta_2)=0.2\qquad m_h(\theta_2\cup\theta_5\cup\theta_6|\theta_2)=0.1$$
\end{itemize}

\subsection{Another example with $\Theta=\{\theta_1,\ldots,\theta_6\}$}

LetÕs change the previous counter-example and use now the following mass matrix:

\begin{center}
\begin{tabular}{|l|c|c|c|c|c|}
\hline
                   &  $\theta_1$ & $\theta_2$ & $\theta_3$ & $\theta_4\cup\theta_5$ & $\theta_5\cup\theta_6$\\
\hline
$m_1(.)$  & $1$ &  $0$ &  $0$ &  $0$ &  $0$\\
$m_{\theta_2}(.)$  & $0$ &  $1$ &  $0$ &  $0$ &  $0$\\
\hline
\end{tabular}
\end{center}

\begin{itemize}
\item Using the Dempster's rule of conditioning, one gets: $m(. | \theta_2) = 0/0$ for all the masses.
\item Using the DSm classical rule, one gets: $m(\theta_1\cap\theta_2 | \theta_2) = 1$, and others 0.
\item If now, one finds out that all intersections are empty (we adopt the Shafer's model), then using the DSm hybrid rule, one gets: $m_h(\theta_1\cup\theta_2|\theta_2)=1$, and others 0.
\end{itemize}

\subsection{Generalization}

Let $\Theta= \{\theta_1,\theta_2,\ldots,\theta_n\}$, where $n \geq 2$, and two basic belief functions/masses $m_1(.)$ and $m_2(.)$ such that there exist $1\leq ( i\neq j) \leq n$, where $m_1(\theta_i)=m_2(\theta_j)=1$, and $0$ otherwise. Then the Dempster's rule of conditioning can not be applied because one gets division by zero. 

\subsection{Example with $\Theta=\{\theta_1,\theta_2,\theta_3,\theta_4\}$ and ignorance}

Let's consider $\Theta=\{\theta_1,\theta_2,\theta_1,\theta_2\}$, one expert and a certain ignorant body of evidence over $\theta_3\cup \theta_4$, with the mass matrix:
\begin{center}
\begin{tabular}{|l|c|c|c|}
\hline
                   &  $\theta_1$ & $\theta_2$ & $\theta_3\cup\theta_4$\\
\hline
$m_1(.)$  & $0.3$ &  $0.7$ &  $0$ \\
$m_{\theta_3\cup\theta_4}(.)$  & $0$ &  $0$ &  $1$\\
\hline
\end{tabular}
\end{center}

\begin{itemize}
\item Using the Dempster's rule of conditioning, one gets $0/0$ for all masses $m(.|\theta_3\cup\theta_4)$.
\item 
Using the classical DSm rule, one gets: $m(\theta_1\cap(\theta_3\cup\theta_4) | \theta_3\cup\theta_4) = 0.3$, $m(\theta_2\cap(\theta_3\cup\theta_4) | \theta_3\cup\theta_4) = 0.7$ and others $0$.
\item
If now one finds out that all intersections are empty (Shafer's model), using the DSm hybrid rule, one gets $m(\theta_1\cup \theta_3\cup\theta_4 | \theta_3\cup\theta_4) = 0.3$, $m(\theta_2\cup\theta_3\cup\theta_4 | \theta_3\cup\theta_4) = 0.7$ and others $0$.
\end{itemize}




\subsection{Generalization}
%
Let $\Theta= \{\theta_1,\theta_2,\ldots,\theta_n,\theta_{n+1},\ldots,\theta_{n+m}\}$, for $n \geq 2$ and $m \geq 2$. Let consider the mass $m_1(.)$, which is a row of its values assigned for $\theta_1,\theta_2,\ldots,\theta_n$, and some unions among the elements $\theta_{n+1}$, \ldots, $\theta_{n+m}$ such that all unions are disjoint with each other. If the second mass $m_A(.)$ is a conditional mass, where $A$ belongs to $\{\theta_1,\theta_2,\ldots,\theta_n\}$ or unions among $\theta_{n+1}$, \ldots, $\theta_{n+m}$, such that $m_1(A)=0$, then the Dempster's rule of conditioning can not be applied because on get division by zero, which is undefined.  [We did not consider any intersection of $\theta_i$ because the Dempster's rule of conditioning doesn't accept paradoxes]. But the DSm rule of conditioning does work here as well.

\subsection{Example with a paradoxical source}

A counter-example with a paradox (intersection) over a non-refinable frame, where the Dempster's rule of conditioning can not be applied because the Dempster-Shafer theory does not accept paradoxist/conflicting information between elementary elements $\theta_i$ of the frame $\Theta$:\\

LetÕs consider the frame of discernment $\Theta= \{\theta_1,\theta_2\}$, one expert and a certain body of evidence over $\theta_2$, with the mass matrix:

\begin{center}
\begin{tabular}{|l|c|c|c|c|}
\hline
                   &  $\theta_1$ & $\theta_2$ & $\theta_1\cap\theta_2$ & $\theta_1\cup\theta_2$\\
\hline
$m_1(.)$  & $0.2$ &  $0.1$ &  $0.4$ &  $0.3$\\
$m_{\theta_2}(.)$  & $0$ &  $1$ &  $0$ &  $0$\\
\hline
\end{tabular}
\end{center}

Using the DSm rule of conditioning, one gets
$$m(\theta_1 | \theta_2) = 0\qquad m(\theta_2 | \theta_2) = 0.1+0.3 = 0.4\qquad m(\theta_1 \cap \theta_2 | \theta_2) = 0.2+0.4 = 0.6\qquad m(\theta_1 \cup \theta_2 | \theta_2) = 0$$
\noindent and the sum of fusion results is equal to 1.\\

Suppose now one finds out that all intersections are empty. 
Using the DSm hybrid rule when $\theta_1 \cap \theta_2=\emptyset$, one has:
%
\begin{align*}
m_h(\theta_1\cap \theta_2|\theta_2) &=0\\
m_h(\theta_1|\theta_2) &= m(\theta_1|\theta_2) + [m_1(\theta_1)m_2(\theta_1\cap \theta_2) + m_2(\theta_1)m_1(\theta_1\cap\theta_2)] = 0\\
m_h(\theta_2|\theta_2) & = m(\theta_2|\theta_2) + [m_1(\theta_2)m_2(\theta_1\cap\theta_2) + m_2(\theta_2)m_1(\theta_1\cap\theta_2)] = 0.4+ 0.1(0)+1(0.4) = 0.8\\
m_h(\theta_1\cup \theta_2|\theta_2) &= m(\theta_1\cup\theta_2|\theta_2) + [m_1(\theta_1)m_2(\theta_2) + m_2(\theta_1)m_1(\theta_2)] \\
& \qquad + [m_1(\theta_1\cap\theta_2)m_2(\theta_1\cup\theta_2) +  m_2(\theta_1\cap\theta_2)m_1(\theta_1\cup\theta_2)] + [m_1(\theta_1\cap\theta_2)m_2(\theta_1\cap\theta_2)]\\
& = 0 + [0.2(1) + 0(0.1)] + [0.4(0) + 0(0.3)] + [0.4(0)] \\
&= 0.2 + [0] + [0] + [0] = 0.2
\end{align*}

\section{ Conclusion}

Several infinite classes of counter-examples to Dempster's rule of combining have been presented in this paper for didactic purpose to show the limitations of this rule in the DST framework. These infinite classes of fusion problems bring the necessity of generalizing the DST to a more flexible theory which permits the combination of any kind of sources of information with any degree of conflict and working on any frame with exclusive or non-exclusive elements. The DSmT with the DSm hybrid rule of combining proposes a new issue to satisfy these requirements based on a new  solid mathematical framework.

\end{document}